\documentclass[12pt,twoside]{amsart}
\usepackage{amssymb}

\newtheorem{theorem}{Theorem}[section]
\newtheorem{lemma}[theorem]{Lemma}
\newtheorem{proposition}[theorem]{Proposition}
\newtheorem{corollary}[theorem]{Corollary}

\theoremstyle{definition}
\newtheorem{definition}[theorem]{Definition} % \theoremstyle{remark}
\newtheorem{remark}[theorem]{Remark}
\newtheorem{example}[theorem]{Example}

\newtheorem{construction}[theorem]{Construction}

\newcommand{\cO}{{\mathcal O}}
\newcommand{\cI}{{\mathcal I}}

\newcommand{\fa}{{\mathfrak a}}
\newcommand{\fb}{{\mathfrak b}}
\newcommand{\fc}{{\mathfrak c}}
\newcommand{\aCM}{arithmetically Cohen-Macaulay}
\newcommand{\lCM}{locally Cohen-Macaulay}

\newcommand{\HR}{Hartshorne-Rao module}
\newcommand{\ga}{\gamma}
\newcommand{\be}{\beta}
\newcommand{\al}{\alpha}

\newcommand{\gin}{\operatorname{gin}}
\newcommand{\iin}{\operatorname{in}}

\newcommand{\HH}{H_{\mathfrak m}}
\newcommand{\im}{\operatorname{im}}
\newcommand{\coker}{\operatorname{coker}}

\newcommand{\eeps}{\varepsilon}

\newcommand{\rank}{\operatorname{rank}_K}
\newcommand{\mif}{\mbox{if} ~}
\newcommand{\s}{\; | \;}

\newcommand{\TR}{\operatorname{Tor}^R}

\newcommand {\ZZ}{\mathbb{Z}}

\newcommand {\PP}{\mathbb{P}}

\begin{document}

\title[Non-degenerate curves with maximal Hartshorne-Rao module]
{Non-degenerate curves with maximal
Hartshorne-Rao module}

\author[Uwe Nagel]{Uwe Nagel$^*$}

% \address{Department of Mathematics,
% University of Kentucky, Lexington, KY 40506-0027, USA}
% \email{uwenagel@ms.uky.edu}

% \date{\today}

\begin{abstract} Extending results for space curves we 
establish bounds for the cohomology of a non-degenerate 
curve in projective $n$-space. As a consequence, for any given 
$n$ we determine all possible pairs $(d, g)$ where $d$ is the  
degree and $g$ is the  (arithmetic) genus of the curve. 
Furthermore, we show that curves attaining our bounds always 
exist and describe properties of these extremal curves. In
particular, we determine the \HR, the generic initial
ideal and the graded Betti numbers of an extremal curve.
\end{abstract}

% \subjclass{Primary 14H50; Secondary ???}
\thanks{$^*$ Department of Mathematics,
 University of Kentucky, Lexington, KY 40506-0027, USA \\
{\it e-mail}: uwenagel@ms.uky.edu}

%%%%%%%%%%%%%%%%%%%%%%%%%%%%%%%%%%%%%
\maketitle
\tableofcontents
%%%%%%%%%%%%%%%%%%%%%%%%%%%%%%%%%%%%%%%%%%%%

\hfill {\it Dedicated to Silvio Greco on the occasion of his 60th
birthday} \\[1ex]

\section{Introduction} \label{section-intro}

Let $\PP^n$ denote the $n$-dimensional projective space over an
algebraically closed field of characteristic zero. By a curve $C
\subset \PP^n$ we will always understand a closed subscheme which
is locally Cohen-Macaulay of pure dimension $1$, thus in
particular without isolated or embedded points. The curve is
called non-degenerate if it is not contained in a hyperplane of
$\PP^n$. The most important invariants of the curve $C$ are its
degree and its (arithmetic) genus. However, a more precise
description of $C$ requires some knowledge of its Hartshorne-Rao
module. Its importance is particularly highlighted in Liaison
theory (cf.\ \cite{Migliore-book}) and in the classification theory of curves
(cf.\ \cite{MDP-book}).

Every graded module of finite length is (up to degree shift) the
Hartshorne-Rao module of a curve (cf.\ \cite{MNP-constr-coho}). However, fixing
the degree $d$ and
the (arithmetic) genus $g$ puts   restrictions on the possible Hartshorne-Rao
modules. In this paper we mainly investigate the restrictions on  the ``size'' of the
Hartshorne-Rao module. More precisely, we consider the following basic problems.
\medskip

\noindent {\bf Problems.} Let $C \subset \PP^n$ be a
non-degenerate curve of degree $d$ and (arithmetic) genus $g$.

(1) Which pairs $(d, g)$ can occur?

(2) Find (optimal) upper estimates for the dimension $h^1(\cI_C(j))$,
$j \in \ZZ$.

(3) Are there curves attaining the bounds in (2) for all $j \in
\ZZ$?
\medskip

The answers to these problems are completely known for space curves, i.e.\ when
$n = 3$.
The goal of this paper is to answer these questions for arbitrary $n
\geq 3$  and to describe properties of the extremal curves occurring in Problem
(3). Hereby we are using tools from cohomology theory, homological algebra
 as well as the theory of Gr\"obner bases. On a technical level
 one  problem is to replace codimension two oriented methods by
 those which work in any codimension. Some of the results in this note
 have been announced in \cite{N-survey}.

The problems above, in particular Problem (1), are very classical for irreducible,
reduced curves. However, we consider them for more general
curves because even if one is primarily interested
in integral curves one often is led to study more generally \lCM\ curves.
This is evident, for example,  in Liaison theory or in the study of families of
curves. For the latter the Hilbert scheme $H^n_{d, g}$ of locally
Cohen-Macaulay curves in
$\PP^n$ of degree $d$ and genus $g$ is the right framework.  Of course, the first
problem above just asks: When does the Hilbert scheme  $H^n_{d, g}$
contain a non-degenerate curve?
\smallskip

Let us describe the organization of the paper. Bounds for the
first and second cohomology of the ideal sheaf of a curve are established in
Section
2. These results lead immediately to a bound on the genus depending on the degree.
In Section 3 we show that this is the only restriction on the genus by
constructing  suitable curves. A refinement of the construction allows to
construct curves attaining the bound on the Hilbert function of the \HR\ in
every degree. This is done in Section 4 and the resulting curves
are called extremal curves.

The remaining parts of this note are
devoted to study these curves. It turns out that the situation for space
curves is more rigid and simpler than for curves of codimension $\geq 3$.
In Section 5 we describe the
structure of the \HR\ of an extremal curve. A key
is the determination of the generic initial ideal of an extremal
curve. Furthermore, we show that extremal curves also have maximal
second cohomology and (mostly) contain a planar  subcurve of degree $d-1$.
The converse to these two results ist only true if $n=3$.

In Section 6 we compute the graded Betti
numbers of an extremal curve. As a consequence we get that its homogeneous
ideal  is componentwise linear if the degree of the curve is at least $5$.

The results in \cite{MDP-Hilb-scheme} strongly suggest that it is
worth to investigate families of extremal curves using the results of Section 5 and 6.
 We wish to pursue  this in a future paper.

%%%%%%%%%%%%%%%%%%%%%%%%%%%%%%%%%%%%%%%%%%%%%%%%%%%%%%%%%%%%%%%%%%%%%%%%5

\section{Bounds for the cohomology} \label{section-bounds}

In this section we establish estimates for the first and second
cohomology of a curve $C$. Note that in \cite{CGN1}, Theorem 2.1 such bounds
have been achieved under the additional assumption that the general
hyperplane section of $C$ in non-degenerate.
\smallskip

If a curve of degree 2 is not arithmetically Cohen-Macaulay then
it is a pair of two skew lines or a double line. Thanks to the
results in \cite{NNS1} we have a complete description of all
double lines. Thus, in this note we often will focus on curves of
degree $d \geq 3$.

\begin{proposition} \label{prop-bounds}
Let $C \subset \PP^n$ be a non-degenerate curve of degree $d \geq 3$
and (arithmetic) genus $g$. Then we have:
\begin{itemize}
\item[(a)]
$$
h^1(I_C(j)) \leq \rho_n^{ex}(j) \quad \mbox{for all} \; j \in \ZZ
$$
\smallskip

where $\rho_n^{ex}: \ZZ \to \ZZ$ is the function defined by
$$
\rho_n^{ex}(j) =
\left \{
\begin{array}{*3l} 0&{\rm if}& j \leq - {d-2 \choose
2} + g\\ \\
{d-2 \choose 2} - g + j& {\rm if} & -
{d- 2 \choose 2} + g \leq j \leq 0\\ \\
{d-2 \choose 2} - g - (n-3) &  {\rm if} & 1 \leq j \leq d-2\\ \\
{d-1 \choose 2} - g - (n-3) - j& {\rm if}&
d-2 \leq j \leq {d -1 \choose 2} - g - (n-3) \\ \\
0&{\rm if}& {d-1 \choose 2} - g - (n-3) \leq j \end{array}
\right.
$$
\item[(b)]
$$
h^2(\cI_C (j)) \leq \mu_n^{ex}(j) \quad \mbox{for all} \, j \in \ZZ
$$
\smallskip

where $\mu_n^{ex}: \ZZ \to \ZZ$ is the function defined by
$$
\mu_n^{ex}(j) =  \left \{
\begin{array}{ll}
0 & \mif j \geq d-3 \\ \\
\binom{d-2-j}{2} & \mif 0 \leq j \leq d-2 \\ \\
\binom{d-2}{2} - (d-1) j - 1  & \mif g - \binom{d-2}{2} \leq j
\leq -1 \\ \\
g - 1 - d j & \mif j \leq g - \binom{d-2}{2} - 1.
\end{array} \right.
$$
\end{itemize}
\end{proposition}

The strategy of proof is similar to the one for Theorem 2.1 in
\cite{CGN1}. Thus, we will focus on the differences and refer to
\cite{CGN1} for more details.

Throughout the paper we will denote by $R := K[x_0,\ldots,x_n]$
the homogeneous coordinate ring of $\PP^n$.
Let $\Gamma := C \cap H$ be the general hyperplane section of $C$.
While $C \subset \PP^n$ is non-degenerate $\Gamma \subset H \cong
\PP^{n-1}$ might be degenerate. The first step consists in estimating
 the failure of non-degeneracy of $\Gamma$. It strengthens
\cite{BCG-lin}, Proposition 4.2.

\begin{lemma} \label{lem-restr}
Let $C \subset \PP^n, n \geq 3,$ be a non-degenerate curve of degree $d
\geq 3$.
Then its general hyperplane section $\Gamma$ is not collinear.
\end{lemma}

\begin{proof}
Let $\ell \in R$ be the linear form defining the general
hyperplane $H$. Put $m := h^0(\cI_{\Gamma}(1))$. If $m = 0$ there
is nothing to show. Assume $m \geq 1$.
Let $f_1, \dots , f_m$ be linearly independent linear
forms in the homogeneous ideal $I_{\Gamma}$ of $\Gamma \subset H$. Then
$T := \overline R/(\ell, f_1 ,\dots,f_m)$ is the homogeneous coordinate
ring of the linear span $L$ of $\Gamma$. Let
$\fc := I_{\Gamma} T$ be the ideal of $\Gamma \subset
L$. Using the so-called Socle lemma of Huneke and Ulrich (cf.\
\cite{HU-hyper}) it is shown in Steps 1 and 3 of the proof of
\cite{CGN1}, Theorem 1.4 that the ideal $\fc$ has codimension
$\leq n-1-m$ and contains at least $n-1-m$ linearly independent
quadrics. It follows that $m \leq n-2$. We have to exclude the
possibility $m = n-2$. Indeed, if $m = n-2$ then we obtain
$\deg C = \deg \Gamma = \deg \fc \leq 2$. This contradiction
completes the proof.
\end{proof}

\begin{remark}
The assumption on the characteristic of the ground field is used
for the application of the Socle lemma. In fact, the statement
above is not true over fields of positive characteristic. There are
multiple lines in $\PP^3$ whose general hyperplane section is
collinear (cf.\ \cite{Hartshorne-genus-1996}).
\end{remark}

The second step in the proof of Proposition \ref{prop-bounds} is
mainly an application of
the methods developed in \cite{B}, \cite{N} and \cite{BN}. For details we refer
to the proof of \cite{CGN1}, Theorem 2.1.

\begin{proof}[Proof of Proposition \ref{prop-bounds}]
We proceed in several steps.
\smallskip

\noindent \underline{Notation}. We use $h_C$ and $p_C$ (resp.\
$h_{\Gamma}$ and $p_{\Gamma}$) to denote the
Hilbert function and the Hilbert
polynomial of $C$ (resp.\ of $\Gamma$).
\medskip

\noindent \underline{Step 1}. According to Lemma \ref{lem-restr}
we have
$h_{\Gamma}(1) \geq 3$, thus $h^1(\cI_{\Gamma}(1))
\leq d - 3$.
Therefore Lemma 4.3 of \cite{N} provides for $j \geq 1$:
$$
h^1(\cI_{\Gamma}(j)) \leq \max \{ 0, d-2-j \}
$$
and
$$
h_{\Gamma}(j) \geq \min \{d, 2 + j \}.
$$
\medskip

\noindent \underline{Step 2}. The long exact cohomology sequence
induced by the exact sequence
$$
0 \rightarrow \cI_C(j - 1) \rightarrow \cI_C(j)
\rightarrow \cI_{\Gamma}(j) \rightarrow 0.
$$
ends as follows:
$$
\cdots \rightarrow H^1(\cI_{\Gamma}(j))
\rightarrow H^2(\cI_C(j - 1)) \rightarrow H^2(\cI_C(j))
\rightarrow 0
$$
Since $H^2(I_C(j)) = 0$ for $j \gg 0$, we get by Step
1 if
$j \geq 0$:
\begin{eqnarray*}
h^2(\cI_C(j)) &\leq &\sum_{t\geq
j+1}h^1(\cI_{\Gamma}(t))\\ &\leq& \sum_{t \geq
j+1} \max \{0, d - 2 - t \} \\ &=& {d - 2 - j \choose 2}\\
\end{eqnarray*}
proving claim (b) for non-negative degrees.
\medskip

\noindent\underline {Step 3}.
Consider now the following version of the Riemann-Roch
theorem:

\begin{equation*}
h_C(j) - p_C(j) = - h^1(I_C(j)) + h^2(I_C(j)) \quad (j
\in \ZZ)
\end{equation*}
where $p_C(j) = d j - g + 1$.

Using Step 2 and $h_C (1) = n + 1$ we obtain for $j \geq 1$:
\begin{eqnarray*}
h^1(\cI_C(j)) &=& p_C(j) - h_C(j) + h^2(\cI_C(j))\\
&\leq& d j - g + 1 - \left [ n+1 + \sum_{t=2}^j
h_{\Gamma}(t) \right ] + \\ & &+ \sum_{t\geq j+1} \max \{0,
d - 2 - t \} \quad {\rm (by\ Step\ 2)}\\ \\
&\leq& d-g-n + \sum_{t\geq 2} \max\{0,d-2-t \} \quad {\rm (by\
Step\ 1)}\\ \\
&=&{d - 2\choose 2} - g + 3 - n.\\
\end{eqnarray*}
For $j = 0$ we get
$$
h^1(\cI_C) \leq - g + {d - 2\choose 2}.
$$
This proves claim (a) if $0 \leq j \leq d-2$. For other values of
$j$ the result follows now by the general methods developed in
\cite{BN} and \cite{N} (cf.\ Steps 4 and 5 in the proof of
\cite{CGN1}, Theorem 2.1).
\medskip

\noindent\underline {Step 4}. For $j < 0$ the Riemann-Roch theorem
implies
$$
h^2(\cI_C (j)) = h^1(\cI_C(j)) - [ d j - g + 1 ].
$$
Thus, in this case claim (b) follows by (a). Together with Step 2 this
completes the proof.
\end{proof}

Observe that part (a) of Proposition \ref{prop-bounds} generalizes the estimates
 of Martin-Deschamps and Perrin in \cite{MDP-bounds} from $n = 3$ to $n \geq 3$
 while part (b) extends \cite{CGN4}, Lemma 4.8.
\smallskip

In Proposition \ref{prop-bounds} we have defined the function $\rho_n^{ex}$ for
curves of degree $d \geq 3$. We want to extend it to curves of degree $2$.

\begin{definition} \label{def-gen-ex-fkt} If  $d = 2$ we set
$$
\rho_n^{ex}(j) :=
\left \{
\begin{array}{*3l} 0&{\rm if}& j \leq  g\\ \\
  - g + j& {\rm if} &  g \leq j \leq 0\\ \\
 - g - (n - 3) - j& {\rm if}&
1 \leq j \leq  - g - (n - 3) \\ \\
0&{\rm if}&  - g - (n-3) \leq j \end{array}
\right.
$$
\end{definition}

Using this notation we can generalize Proposition \ref{prop-bounds}

\begin{theorem} \label{thm-bounds}
Let $C \subset \PP^n,\ n \geq 3,$ be a non-degenerate curve of degree $d$
and (arithmetic) genus $g$. Then we have:
$$
h^1(I_C(j)) \leq \rho_n^{ex}(j) \quad \mbox{for all} \; j \in \ZZ.
$$
\end{theorem}

\begin{proof}
Since the curve $C$ is non-degenerate its degree must be at least two.
 Thus, by virtue of Proposition \ref{prop-bounds}
it remains to consider the case $d=2$. If $C$ is \aCM\ then it is
a planar conic contradicting our assumption. Hence, $C$ is not \aCM, thus it
 must be a pair of two skew lines or a double line. Therefore,  the claim follows by
\cite{NNS1}, Corollary 3.2.
\end{proof}

The objective of the next two sections is to show that the bounds
in the theorem above cannot be improved.

%%%%%%%%%%%%%%%%%%%%%%%%%%%%%%%%%%%%%%%%%%%%%%%%%%%%%%%%%%%%%%%%%%%%%%%%%

\section{The possible genera of curves} \label{section-genera}

Theorem \ref{thm-bounds} immediately implies.

\begin{corollary}
Let $C \subset \PP^n$ be a non-degenerate curve of degree $d \geq
3$. Then the arithmetic genus $g$ of $C$ satisfies
$$
g \leq \binom{d-2}{2} - (n-3).
$$
\end{corollary}

Here we want to show that this is the only restriction. We also include curves of
degree $2$. Notice, that in this case \cite{NNS1}, Corollary 2.10 provides $ g \leq
2-n$.

\begin{theorem} \label{thm-genus-range}
A non-degenerate curve $C \subset \PP^n, n \geq 3,$ of degree
$d \geq 2$ can have every \\[3pt]
genus $g \leq g_{max} := \left \{ \begin{array}{cl}
\binom{d-2}{2} - (n-3) & \mif d \geq 3 \\[1ex]
2 -n & \mif d =2.
\end{array} \right. $
\end{theorem}

We will show this result by constructing suitable curves. We use
ideas of the proof of Proposition 1.11 in \cite{CGN1} where the
statement has been shown in case $n=4$. For $n=3$ the result has
been proved independently in \cite{Sauer}, \cite{Okonek-Sp} and
\cite{Hartshorne-genus-1996}.

\begin{construction} \label{constr}
Let $D \subset \PP^n, n \geq 3,$ be a planar curve of degree $d-1$ which is
supported on the line $L$. Let $a \geq 0$ be an integer.

We may assume that the homogeneous ideals of $L$ and $D$ in $R =
K[x_0,\ldots,x_n]$ are $I_L = (x_2,\ldots,x_n)$ and
$I_D = (x_2^{d-1}, x_3,\ldots,x_n)$, respectively. Let $S = K[x_0, x_1] \cong
R/I_L$ be the homogeneous coordinate ring of $L$. Observe that
there is an isomorphism of graded $S$-modules
$$
I_D/I_L I_D \cong S^{n-2}(-1) \oplus S(-d+1).
$$
Now we choose homogeneous elements $f_1,\ldots,f_{n-2}, f \in S$
such that $f_1,\ldots,f_{n-2}$ have degree $a+n-3$ and $f$ has
degree $d+a+n-5$. Notice that our assumptions imply $a+n-3 \geq
0$ and $d+a+n-5 \geq 0$. Using these forms we define the following
homomorphism
$$
\psi: S^{n-2}(-1) \oplus S(-d+1)
\stackrel{(f_1,\ldots,f_{n-2},f)}{\longrightarrow} S(a+n-4).
$$
Let $\al: I_D \to S(a+n-4)$ be the composition
$$
I_D \to I_D/I_L I_D \stackrel{\cong}{\longrightarrow} S^{n-2}(-1) \oplus
S(-d+1) \stackrel{\psi}{\longrightarrow} S(a+n-4)
$$
where the first map is the canonical epimorphism. Since $D$ and
$L$ are \aCM\ $\ker \al \subset I_D$ is a homogeneous saturated
ideal defining a $1$-dimensional subscheme $C \subset \PP^n$.
Moreover, we have an exact sequence of graded $R$-modules
$$
0 \to I_C \to I_D \stackrel{\al}{\longrightarrow} R/I_L (a+n-4)
\to \coker \al \to 0. \leqno(*)
$$
\end{construction}
\medskip

Keeping the notation and assumptions of the construction we get.

\begin{lemma} \label{lem-constr-sequence}
Assume that $\coker \al$ has finite length, $d \geq 3$ and that
$f_1,\ldots,f_{n-2}$ are linearly independent. Then there is a
non-degenerate curve $C \subset \PP^n$ fitting into the exact sequence
$$
0 \to  \cO_L(a+n-4) \to \cO_C \to \cO_D \to 0.
$$
\end{lemma}

\begin{proof}
Since $\coker \al$ has finite length, sheafification of the  exact
sequence $(*)$ above provides the exact sequences
$$
0 \to \cI_C \to \cI_D \to \cO_L(a+n-4) \to 0
$$
and
$$
0 \to  \cO_L(a+n-4) \to \cO_C \to \cO_D \to 0.
$$
Clearly, the support of $C$ is $L$. Furthermore, the second
sequence shows for all closed points $P \in L$ that $\cO_{C, P}$ has
depth $1$. Hence is $C$ is a locally Cohen-Macaulay subscheme of
pure dimension $1$, i.e.\ a curve.

Since $\ker \al$ is the homogeneous ideal of $C$, the curve $C$
is non-degenerate if and only if $[\ker \psi]_1 = 0$. But this
just means that $f_1,\ldots,f_{n-2}$ are linearly independent
because $[S(-d+1)]_1 = 0$.
\end{proof}

Now we can show Theorem \ref{thm-genus-range}  easily.

\begin{proof}[Proof of Theorem \ref{thm-genus-range}]
If $d = 2$ the claim follows by \cite{NNS1}, Theorem 2.6 and Proposition 2.9.

Assume $d \geq 3$. Of course,
we want to apply the lemma above. We have to check that we
can satisfy its assumptions. But $[S]_{a+n-3}$ is a $K$-vector
space of dimension $a+n-2 \geq n-2$. Thus we can find linearly
independent elements $f_1,\ldots,f_{n-2}$ in $[S]_{a+n-3}$.
Moreover, we can certainly arrange that $\coker \psi$ and thus $\coker
\al$ has finite length. For example, choose $f_1 = x_0^{a+n-3}$ and $f =
x_1^{d+a+n-5}$.

Now using the sequence in Lemma \ref{lem-constr-sequence} an easy
computation of Hilbert polynomial
shows that $C$ is a curve of degree $d$ and arithmetic genus
$g_{max}- a$. Since $a$ is any
non-negative integer we are done.
\end{proof}

\begin{remark} An analysis of the construction above shows that in
order to prove Theorem \ref{thm-genus-range} we could even
construct a curve $C$ having the two additional properties:
\begin{itemize}
\item[(a)]  $C$ has generically embedding dimension 2 (in
particular it is a generically complete
intersection);
\item[(b)] $C$ is irreducible, i.e.\ it is not the scheme-theoretic union of any
two closed subschemes of $C$, different from $C$.
\end{itemize}
This follows as in the proof of \cite{CGN1}, Proposition 1.11.
 The reader interested in seeing such a curve described by explicit
equations is referred to Example \ref{ex-equations}.
\end{remark}

%%%%%%%%%%%%%%%%%%%%%%%%%%%%%%%%%%%%%%%%%%%%%%%%%%%%%%%%%%%%%%%%%%%%%%%%%

\section{Construction of extremal curves} \label{section-constr}

Following Martin-Deschamps and Perrin  we define.

\begin{definition} \label{def-extremal}
A non-degenerate curve $C \subset \PP^n$ is
called {\it extremal} if
$$
h^1(\cI_C(j)) = \rho_n^{ex}(j) \quad \mbox{for all} \; j \in \ZZ.
$$
\end{definition}

A priori it is not clear at all that such extremal curves do
exist. In case $ n = 3$ the existence has been established by
Martin-Deschamps and Perrin (\cite{MDP-bounds}). A first generalization of this
result to $n \geq 4$ has been achieved in \cite{CGN1} where it was
assumed that with $C$ also its general hyperplane section is
non-degenerate. Here we want to generalize the result for $n=3$
 to $n \geq 3$ without assuming additional properties. This is taken
care of by refining the
construction of the previous section.

\begin{theorem} \label{thm-exist-extre}
For every pair $(d, g)$ of integers such that $d \geq 2$ and
$g \leq g_{max}$ there is an extremal curve $C \subset \PP^n$.
\end{theorem}

\begin{remark}
Contrary to the situation for space curves, extremal curves of
$\PP^n$ of maximal genus $g = g_{max}$ are not \aCM\ if $n \geq
4$.
\end{remark}

Much of this section is based on the following observation concerning
Construction \ref{constr}.

\begin{lemma} \label{lem-HR-constr}
With the notation and assumptions of Construction \ref{constr}
assume additionally that $\coker \al$ has finite length. Then the
scheme $C$ constructed there is a curve with \HR\
$$
H^1_*(\cI_C) \cong (R/(I_L + (f_1,\ldots,f_{n-2},f) R)) (a+n-4).
$$
\end{lemma}

\begin{proof}
First of all note that
$$
\coker \al = \coker \psi \cong (R/(I_L + (f_1,\ldots,f_{n-2},f) R)) (a+n-4).
$$
The exact sequence $(*)$ in Construction \ref{constr} provides the
exact sequences
$$
0 \to I_C \to I_D \to \im \al \to 0
$$
and
$$
0 \to \im \al \to  R/I_L (a+n-4)
\to \coker \al \to 0.
$$
Since $R/I_L$ is Cohen-Macaulay the second sequence implies
$$
\HH^1(\im \al) \cong \HH^0(\coker \al) \cong \coker \al.
$$
Using the first sequence and that $D$ is \aCM\ we obtain
$$
\HH^2(I_C) \cong \HH^1(\im \al).
$$
Taking into account the isomorphism
$$
H^1_*(\cI_C) \cong \HH^2(I_C)
$$
the claim follows.
\end{proof}

This lemma reduces the proof of Theorem \ref{thm-exist-extre} to a
suitable choice of the polynomials $f_1,\ldots,f_{n-2}, f$.

\begin{proof}[Proof of Theorem \ref{thm-exist-extre}]
If $d = 2$ the claim follows by \cite{NNS1}, Corollary 3.2.

Assume $d \geq 3$.
Let $h \neq 0$ be a homogeneous polynomial of degree $a := g_{max} - g \geq 0$ in
$S = k[x_0, x_1]$. Then we choose $f_1,\ldots,f_{n-2} \in [S]_{a+n-3}$ such that
$$
(f_1,\ldots,f_{n-2}) = h\ (x_0, x_1)^{n-3}.
$$
Finally, we choose $f \in [S]_{d+a+n-5}$ such that $\{f, h\}$ is a
regular sequence unless $a = 0$ and $n = 3$. (If  $a = 0$  there is no condition of $f$.)

Then it is easy to see that we have for all integers $j \in \ZZ$
$$
\rank [R/(I_L + (f_1,\ldots,f_{n-2},f)]_{j+a+n-4} = \rho^{ex}_n
(j).
$$
Therefore Lemma \ref{lem-HR-constr} and Lemma
\ref{lem-constr-sequence} show that $C$ is an extremal curve.
\end{proof}

For each pair $(d, g)$ allowed by Theorem \ref{thm-exist-extre}
we describe one extremal curve by explicit equations.

\begin{example} \label{ex-equations}
In the last proof we specialize the polynomials $f, h$ to
$$
f := x_1^{d+a+n-5} \quad \mbox{and} \quad h := x_0^a
$$
where $a := g_{max} - g$. Then, it is not too difficult to see that $\ker \alpha$ is
the ideal
$$
\begin{array}{rcl}
I &  = & (x_2^{d-1}, x_3,\ldots,x_n) \cdot (x_2,\ldots,x_n) \;  + \\[1ex]
& & \hspace*{2cm}
(x_0^a x_2^{d-1} + x_1^{d+a-2} x_3, x_0 x_i + x_1 x_{i+1} \s 3 \leq i \leq n-1).
\end{array}
$$
Hence $I$ defines an extremal curve $C \subset \PP^n$ of degree $d$ and genus $g =
g_{max} - a$ provided $d \geq 3$. If $d = 2$ and $a > 0$ then the ideal $I$ defines an
extremal curve as well. It has genus $3-n-a$.  Notice that in all cases this
extremal curve  is irreducible and has generically embedding dimension
$2$.
\end{example}

We want to point out that the situation is more rigid for $n = 3$
than for $n \geq 4$. Indeed, it follows by
\cite{Nollet-subex-curves} that a non-degenerate curve $C \subset
\PP^3$ of degree $d \geq 4$  with $h^1(\cI_C(1)) =
\rho^{ex}_3 (1) > 0$ must be an extremal curve. The  example below
shows that the analogous conclusion is not true if $n \geq 4$. It
is  even false if we assume $h^1(\cI_C(j)) =
\rho^{ex}_3 (j)$ for all $j \leq 1$.

\begin{example} \label{ex-non-max-Rao}
We use again the Construction \ref{constr}.
Suppose that $n \geq 4$ and $a > 0$. Define integers $k, \eeps$ by
$$
a =: (k-1) (n-3) + \eeps \quad \mbox{where} \; 0 \leq \eeps \leq n-4.
$$
Put $b =: (n-2) (a+n-3)$ and define the polynomials $f_1,\ldots,f_{n-2}$ by
$$
(f_1,\ldots,f_{n-2}) = (x_0^b, x_0^{i k} x_1^{(n-3-i) k + \eeps} \s 0 \leq i
\leq n-4).
$$
Setting $f := 0$ (which is allowed!) we obtain by Lemma
\ref{lem-constr-sequence}  a non-degenerate curve $C \subset \PP^n$ of
degree $d$ and genus $g_{max} - a$ whose Hartshorne-Rao module is given by
Lemma \ref{lem-HR-constr}. Let $d \geq 4$. Then we get
$$
h^1(\cI_C(j)) = \rho^{ex}_3 (j) \quad \mbox{for all} \; j \leq 1
$$
and
$$
h^1(\cI_C(2)) < h^1(\cI_C(1)) = a = \rho_n^{ex}(1) = \rho_n^{ex}(2).
$$
Hence $C$ is not an extremal curve.
\end{example}

\begin{remark} \label{rem-possible-Rao}
We have seen that the Construction \ref{constr} produces a non-degenerate
curve of degree $d \geq 3$ and genus $g_{max} - a$ for all choices of linearly
independent forms $f_1,\ldots,f_{n-2} \in S$ of degree $a+n-3$ and a form
$f \in S$ of degree $d+a+n-5$ such that $S/(f, f_1,\ldots,f_{n-2})S$ has
finite length. Then Lemma \ref{lem-HR-constr} gives for the resulting curve
$C$ that $H^1_*(\cI_C) \cong (S/(f, f_1,\ldots,f_{n-2})S) (a+n-4)$. Thus, it
is not too difficult to determine the possible Rao functions of the curves
obtained by Construction \ref{constr}. In particular, every such curve $C
\subset \PP^n$ satisfies
$$
h^1(\cI_C(j)) = \rho^{ex}_n (j) \quad \mbox{for all} \; j \leq 1.
$$
\end{remark}
\smallskip

The next remark shows that the situation in case $d=3$ is somewhat particular.

\begin{remark} \label{rem-d3-more-HR}
We will see in the next section that Lemma \ref{lem-HR-constr} describes the \HR\ of
{\it every} extremal curve of degree $d \geq 3$ by its generator and minimal relations
unless $d=3, n\geq 4$ and $a > 0$. Indeed, in this case Lemma \ref{lem-HR-constr}
and Lemma \ref{lem-constr-sequence} show that, for example, also
$$
R/(x_2,\ldots,x_n, x_1^{a+n-3}, x_0^{a+i} x_1^{n-3-i} \s 1 \leq i \leq n-3)
$$
is the \HR\ module of an extremal curve of degree $3$ and genus $g_{max} - a$.
\end{remark} 

From an Liaison-theoretic point of view minimal curves are particular 
interesting. In this respect we get. 

\begin{remark} \label{rem-min-curves} 
Since every curve $C$  has the property that the function
$j \mapsto h^1(\cI_C(j))$ is strictly increasing for $j =
t,\ldots,0$ provided $h^1(\cI_C(t)) > 0$ where $t < 0$, every
extremal curve is minimal in its even Liaison class.
\end{remark} 

In the next section we will show that extremal curves also
have maximal second cohomology.

%%%%%%%%%%%%%%%%%%%%%%%%%%%%%%%%%%%%%%%%%%%%%%%%%%%%%%%%%%%5

\section{Properties of extremal curves} \label{section-propert}

The construction of the previous section does not provide all extremal curves.
In fact, we do not have an explicit description of all such curves.
Nevertheless, we are able to derive properties of extremal curves.
Our results show that the situation for space
curves is more rigid and simpler than for curves of codimension $\geq 3$.
\smallskip

We begin with the computation of the \HR.

\begin{theorem} \label{thm-HR-mod-extr}
Let $C \subset \PP^n$ be a non-degenerate curve of degree $d \geq
3$ and genus $g$. Suppose that at least one of the three conditions
\begin{itemize}
\item[(i)] $d = 3$,
\item[(ii)] $g < g_{max} = \binom{d-2}{2}- (n - 3)$,
\item[(iii)] $n \geq 4$
\end{itemize}
is not satisfied.
Then $C$ is an extremal curve if and only if its
Hartshorne-Rao module is (up to change of coordinates) isomorphic
to $R/(x_2,\ldots,x_n, h \cdot (x_0, x_1)^{n-3}, f) (\binom{d-2}{2} - g -1)$
where  $\deg h = g_{max}  - g, \deg f = \binom{d-1}{2} - g$ and
$\{x_2,\ldots,x_n, f, h\}$ is a regular sequence if $g < g_{max}$ and $h = 1$
if $g = g_{max}$.
\end{theorem}

\begin{remark}
Remark \ref{rem-d3-more-HR} shows that the theorem cannot be extended to extremal
curves of degree $3$ excluded in the statement above.
\end{remark}

The proof of Theorem \ref{thm-HR-mod-extr} requires some preparation.

\begin{lemma} \label{lem-cyclic-HR}
The \HR\ of an extremal curve is a cyclic $R$-module.
\end{lemma}

\begin{proof}
If $d = 2$ the claim is a consequence of \cite{NNS1}, Proposition 3.1.

Now assume $d \geq 3$. Put as before $a := g_{max} - g$. We proceed in
two steps and use ideas of \cite{MDP-Hilb-scheme}.
\smallskip

\noindent\underline {Step 1}. Let $l \in R$ be a general linear
form. We show that $[H^1_*(\cI_C)/l H^1_*(\cI_C)]_j = 0$ for all $j
\geq 1$.

To this end denote by $\Gamma$ the intersection of $C$ and the
hyperplane $H$ defined by $l$. Step~1 in the proof of
Proposition \ref{prop-bounds} shows $H^1_* (\cI_{\Gamma}(j)) = 0$ if
$j \geq d-2$.  Since $H^1_*(\cI_C)/l H^1_*(\cI_C)$
is a submodule of $H^1_*(\cI_{\Gamma})$, the claim follows for
$j \geq d-2$. It remains to consider the integers $j$ with $1
\leq j \leq d-3$.

Since $C$ has maximal cohomology, the inequalities in Steps 2 and 3 of
the proof of Proposition \ref{prop-bounds} must be equalities
if $1 \leq j \leq d-2$. This provides the claim for
$j$ in this range.
\smallskip

\noindent\underline {Step 2}. The Hartshorne-Rao module of $C$
must have the form
$$
H^1_*(\cI_C) \cong (R/\fa)(a+n-4)
$$
where $\fa$ is an ideal containing $n-1$ linearly independent
forms.

Step 1 provides that the minimal generators of $H^1_*(\cI_C)$ (as
$R$-module) have degree $\leq 0$. Thus, the claim is clear if $a+n-3 \leq 1$.

Assume $a+n-3 \geq 2$.
The Hilbert function of the $R$-module $H^1_*(\cI_C)$ implies that it
has exactly one
minimal generator, say $y$, of degree $4-n-a$. Denote by $\fa$ the
annihilator of $y$. Then there is an embedding $(R/\fa)(a+n-4) \to
H^1_*(\cI_C)$. It follows that $\dim_K [R/\fa]_1 \leq h^1
(\cI_C(-a-n+5)) = 2$. Thus,
$\fa$ contains $n-1$ linearly independent forms. Suppose $\fa$ has
a further minimal generator $p$ of degree $e \leq a+n-4$. Then we may
assume that the hypersurface $P$ defined by $p$ intersects the
curve $C$ properly. Hence the exact sequence
$$
0 \to \cI_C(-e) \stackrel{p}{\longrightarrow} \cI_C \to \cI_{C \cap P, P} \to 0
$$
implies the exact cohomology sequence
$$
\ldots \to H^0(\cI_{C \cap P, P}(4-a-n + e)) \to H^1(\cI_C(4-a-n))
\stackrel{p}{\longrightarrow} H^1(\cI_C(4-a-n+e)) \to \ldots
$$
where the multiplication by $p$ is the zero map. But this gives a
contradiction since $4-a-n+e \leq 0$, thus $H^0(\cI_{C \cap P, P}(4-a-n +
e))= 0$, but $h^1 (\cI_C(4-a-n)) = 1$. It follows that the minimal
generators of $\fa$ besides
the $n-1$ linear forms  must have degree $\geq a+n-3$. Comparing vector
space dimensions we see that the embedding $(R/\fa)(4-n-a) \to
H^1_*(\cI_C)$ is bijective in all degrees $\leq 0$. Therefore, $y$
is the only minimal generator of $H^1_*(\cI_C)$ and the embedding
is an isomorphism.
\end{proof}

\begin{remark} \label{rem-hilb-v-stru}
By assumption we know the Hilbert function of the \HR\ of an extremal
curve. Thus, one might wonder if this information combined with the
lemma above suffices to determine the module structure. However, this
is not the case as shown by the following example:

Consider in the polynomial ring $S = K[x, y]$ the ideals
$$
\fa = (x^a \cdot (x, y)^{n-3}, y^b)
$$
and
$$
\fb = (x^a \cdot (x, y)^{n-3}, x^{a-1} y^{b-a+1}, y^{b+1})
$$
where $a \geq 2$ and $b > a + n -3$. It is easy to see that $S/\fa$ and $S/\fb$ have
the same Hilbert function. However Theorem \ref{thm-HR-mod-extr} says that (up to
degree shift) $S/\fa$ is, but  $S/\fb$ is not the \HR\ of an extremal curve.
\end{remark}

In order to overcome the problem pointed out in Remark \ref{rem-hilb-v-stru}
we compute the generic initial
ideal of an extremal curve. Hereby we always are using the reverse
lexicographic order on the polynomial ring $R$. It is defined by
\begin{quote}
$x_0^{a_0} \cdot \ldots \cdot x_n^{a_n} > x_0^{b_0} \cdot \ldots \cdot
x_n^{b_n}$ if the last non-zero coordinate of the vector
$$
\left (a_0 - b_0,\ldots,a_n - b_n, \sum_{i=0}^n (b_i - a_i) \right )
$$
is negative.
\end{quote}
We will use various results of the theory of generic initial ideals for
which we refer to \cite{Eisenbud-book} and \cite{Green-gin}. In particular, the
generic initial ideals are stable because the ground field $K$ has characteristic
zero by assumption.

\begin{proposition} \label{prop-gin-ex}
Let $C \subset \PP^n$ be an extremal curve of degree $d \geq 3$. Then its
generic initial ideal is
$$
\gin I_C = (x_0,\ldots,x_{n-4}) \cdot (x_0,\ldots,x_{n-1}) + (x_{n-3}^2, x_{n-3}
x_{n-2}, x_{n-2}^d, x_{n-2}^{d-1} x_{n-1}^a)
$$
where $a := g_{max} - g$ and the first summand is defined to be zero if $n=3$, unless
$d=3, a \geq 1$ and $n \geq 4$. In the latter case the generic initial ideal of $C$
is either the ideal above or
$$
(x_0,\ldots,x_{n-5}) \cdot (x_0,\ldots,x_{n-1}) + x_{n-4} \cdot (x_0,\ldots,x_{n-2}) +
(x_{n-3}^2, x_{n-3} x_{n-2}, x_{n-2}^2, x_{n-4} x_{n-1}^{a+1}).
$$
\end{proposition}

\begin{proof}
Denote by $\Gamma \subset H$ the general hyperplane section of the curve $C$
where $H$ is defined by the linear form $l$. We may assume
that $C$ and $\Gamma$ are given in generic coordinates. Thus, their initial ideals
agree with their generic initial ideals.

We begin by computing $\gin I_{\Gamma} \subset K[x_0,\ldots,x_{n-1}] \cong R/x_n R$. By the proof of
Proposition \ref{prop-bounds} we know that the Hilbert function of $\Gamma$ is
$$
h_{\Gamma} (j) = \min \{j+2, d\} \quad \mif j \geq 1.
$$
Hence $I_{\Gamma}$ contains exactly $n-3$ linear independent linear forms. By
stability of $\gin I_{\Gamma}$ we obtain
$$
x_0,\ldots,x_{n-4} \in \gin I_{\Gamma}.
$$
Since $I_{\Gamma}$ is saturated $x_{n-1}$ is not a divisor of any minimal generator
of $\gin I_{\Gamma}$.

Suppose that $d \geq 4$. Then $I_{\Gamma}$ contains all except four quadrics
of $R/x_n R$. Using stability again we obtain
$$
x_{n-2}^2, x_{n-3} x_{n-1}, x_{n-2} x_{n-1}, x_{n-1}^2 \notin \gin I_{\Gamma}.
$$
If $d =3$ then we get
$$
(x_0,\ldots,x_{n-4}, x_{n-3}^2, x_{n-3} x_{n-2}, x_{n-2}^2) \subset \gin I_{\Gamma}.
$$
Since both ideals have the same Hilbert function they must be equal.

Turning again to the case $d \geq 4$ the Hilbert functions of $\Gamma$ implies that
$I_{\Gamma}$ has a minimal generator of degree $d-1$. Its initial monomial must be
$x_{n-2}^{d-1}$. Comparing Hilbert functions we obtain
$$
(x_0,\ldots,x_{n-4}, x_{n-3}^2, x_{n-3} x_{n-2}, x_{n-2}^{d-1}) = \gin I_{\Gamma}
$$
which is also true if $d =3$ as shown above.

In the next step we compute the quadrics contained in $\gin I_{C}$. Put
$\bar{I_C} := I_C/l I_C$. By Step 1 in the proof of Lemma \ref{lem-cyclic-HR} we
have for $j \geq 1$ an exact sequence
$$
0 \to [I_{\Gamma}/\bar{I_C}]_j \to H^1(\cI_C(j-1)) \stackrel{l}{\longrightarrow}
H^1(\cI_C(j)) \to 0.
$$
It implies
$$
\dim_K [I_{\Gamma}/\bar{I_C}]_j = \left \{ \begin{array}{ll}
n-3 & \mif j = 1 \\
0 & \mif 2 \leq j \leq d-2 \\
1 & \mif d-1 \leq j \leq a+d-2 \\
0 & \mif a+d-1 \leq j.
\end{array} \right. \leqno(+)
$$
Thus, if $a = 0$ then $I_{\Gamma}$ and $\bar{I_C}$ agree in all degrees $j \geq 2$. Our
claim follows by the knowledge of $\gin I_{\Gamma}$.

If $d \geq 4$ and $a\geq 1$ then all the quadrics in
$I_{\Gamma}$ are just restrictions of quadrics in $I_C$. It follows
$$
[\gin I_C]_2 = [\iin I_C]_2 = [((x_0,\ldots,x_{n-4}) \cdot (x_0,\ldots,x_{n-1}) + (x_{n-3}^2, x_{n-3}
x_{n-2})]_2.
$$
Using the dimension count $(+)$ we see that $\gin I_C$ has no further minimal
generators in degrees $\leq d-1$, but one in degree $d$ which is in the ideal
$x_{n-2}^{d-1} R$. By stability we get $x_{n-2}^d \in \gin I_C$. Comparing Hilbert
functions we see that $\gin I_C$ has one further minimal generator of degree $d+a-1$
contained in $x_{n-2}^{d-1} R$. Thus, stability provides
$$
(x_0,\ldots,x_{n-4}) \cdot (x_0,\ldots,x_{n-1}) + (x_{n-3}^2, x_{n-3}
x_{n-2}, x_{n-2}^d, x_{n-2}^{d-1} x_{n-1}^a) \subset \gin I_C.
$$
Since both ideals have the same Hilbert functions we get equality as claimed.

Finally, consider the case $d = 3$ and $a \geq 1$.   Then $I_C$ contains one
quadric less than $I_{\Gamma}$. For reasons of stability it follows that
 $x_{n-2}^2$ or  $x_{n-4} x_{n-1}$, provided $n \geq 4$, is the quadric in
 $\gin I_{\Gamma}$ which is not in $\gin I_C$. In the first case we get as in
 case $d \geq 4$ that
 $$
\gin I_C = (x_0,\ldots,x_{n-4}) \cdot (x_0,\ldots,x_{n-1}) + (x_{n-3}^2, x_{n-3}
x_{n-2}, x_{n-2}^3, x_{n-2}^{2} x_{n-1}^a).
 $$
In the second case we conclude by $(+)$ that $x_{n-4} x_{n-1}^{a+1} \in \gin I_C$.
Comparing Hilbert functions we obtain
$$
\begin{array}{lcl}
\gin I_C & = & (x_0,\ldots,x_{n-5}) \cdot (x_0,\ldots,x_{n-1}) +
x_{n-4} \cdot (x_0,\ldots,x_{n-2}) \\[2pt]
& & + (x_{n-3}^2, x_{n-3} x_{n-2}, x_{n-2}^2, x_{n-4} x_{n-1}^{a+1})
\end{array}
$$
completing the proof.
\end{proof}

We will see in Remark \ref{rem-gins} that both possibilities for the
generic initial ideal of an
extremal curves  of degree $3$ do really occur.
\smallskip

A first consequence of the last result is that extremal curves
also have maximal second cohomology. We use the notation of
Proposition \ref{prop-bounds}.

\begin{corollary} \label{cor-2-extr} Let $C \subset \PP^n$ be an extremal
of degree $d \geq 3$. Then we have
$$
h^2(\cI_C(j)) = \mu^{ex}_n (j) \quad \mbox{for all} \; j \in \ZZ.
$$
\end{corollary}

\begin{proof}
By assumption we know $h^1(\cI_C(j))$ and Proposition
\ref{prop-gin-ex} provides the Hilbert function of $C$. Hence, our
claim is a consequence of the Riemann-Roch formula
$$
h_C(j) - p_C(j) = - h^1(\cI_C(j)) + h^2(\cI_C(j)).
$$
\end{proof}

\begin{remark} \label{rem-2-extr}
(i) An analogous statement is also true for extremal curves of
degree $2$ after modifying the definition of $\mu^{ex}_n$.  In case $d = 2$ we
have to put $\mu^{ex}_n (j) = 0$ if $j \geq 0$. We
leave the details to the reader.

(ii) The last results extends part of \cite{CGN4}, Lemma 4.8 which
also shows that the converse of Corollary \ref{cor-2-extr} is true
provided $n=3$ and $d \geq 5$. However, the curves in Example
\ref{ex-non-max-Rao} show that in case $n \geq 4$ there are curves
of arbitrarily large degree for which the converse of Corollary
\ref{cor-2-extr} is false.
\end{remark}

As the final preparation for the proof of Theorem
\ref{thm-HR-mod-extr} we compute the Betti numbers of the ideals encountered in
Proposition \ref{prop-gin-ex}.

\begin{lemma} \label{lem-res-stabil}
For $1 \leq i \leq n$ define integers
$$
\begin{array}{c}
\al_i  = (n-3) \binom{n}{i} +\binom{n-1}{i} - \binom{n-2}{i+1} \\[1ex]
\al'_i  = \al_i - \binom{n-2}{i-2}, \quad
\be_i  =  \binom{n-2}{i-1}, \quad
\ga_i  =  \binom{n-1}{i-1}.
\end{array}
$$
{\rm (a)} \; Let $a \geq 0, d \geq 3$ be integers.
The ideal
$$
 I = (x_0,\ldots,x_{n-4}) \cdot (x_0,\ldots,x_{n-1}) + (x_{n-3}^2, x_{n-3}
x_{n-2}, x_{n-2}^d, x_{n-2}^{d-1} x_{n-1}^a)
$$
has a minimal free resolution of the form
$$
0 \to R^{n-3} (-n-1) \oplus R(-d-a-n+2) \to F_{n-1} \to \ldots \to F_1 \to I \to 0
$$
where if $a > 0$
$$
F_i = R^{\al_i}(-i-1) \oplus R^{\be_i}(-i-d+1) \oplus R^{\ga_i}(-i-d-a+2) \quad (1
\leq i \leq n-1)
$$
and if $a = 0$
$$
F_i = R^{\al_i}(-i-1) \oplus  R^{\be_i}(-i-d+2) \quad (1
\leq i \leq n-1).
$$

{\rm (b)} \; For $n \geq 4$ and $a > 0$ the minimal free resolution of the ideal
$$
J = (x_0,\ldots,x_{n-5}) \cdot (x_0,\ldots,x_{n-1}) + x_{n-4} \cdot (x_0,\ldots,x_{n-2}) +
(x_{n-3}^2, x_{n-3} x_{n-2}, x_{n-2}^2, x_{n-4} x_{n-1}^{a+1}).
$$
has the shape
$$
0 \to R^{n-4} (-n-1) \oplus R(-a-n-1) \to G_{n-1} \to \ldots \to G_1 \to J \to 0
$$
where
$$
G_i = R^{\al'_i}(-i-1)  \oplus R^{\ga_i}(-i-a-1) \quad (1
\leq i \leq n-1).
$$
\end{lemma}

\begin{proof}
Observe that the two ideals are stable. Thus their minimal free resolution is
described in a theorem of Eliahou and Kervaire \cite{Eliahou-K}. A
straightforward, but lengthy computation provides the graded Betti numbers using
the formula in \cite{Green-gin}, Corollary 1.32 and the identity
$$
\sum_{k = i+1}^{n-3} (k+1) \binom{k}{i-1} = (n-2) \binom{n-2}{i} -
\binom{n-2}{i+1}.
$$
\end{proof}

Combining the information we have obtained we finally can determine the structure of
the \HR\ of extremal curves.

\begin{proof}[Proof of Theorem \ref{thm-HR-mod-extr}]
The other implication being clear it suffices to show that the \HR\ of an extremal
curve $C \subset \PP^n$ (with admissible $d, g, n$) is
$$
H^1_*(\cI_C) \cong R/(x_2,\ldots,x_n, h \cdot (x_0, x_1)^{n-3}, f)
(4-n-a)
$$
where $a = g_{max} - g$.

Using the notation of Step 2 in the proof of Lemma \ref{lem-cyclic-HR} we have
$$
R/\fa \cong S/\fb
$$
where $\fb$ is an ideal in the polynomial ring  $S = K[x, y] \cong R/(x_2,\ldots,x_n)
R$. Thus, it remains to determine the ideal $\fb$. We use the knowledge of the
Hilbert function of $S/\fb$. If $a = 0$ we immediately get $\fb = (x, y)^{n-3}$ as
claimed. Thus, we may assume $a \geq 1$ for the rest of the proof.
Then the Hilbert function of $S/\fb$ still
implies two facts:
\begin{itemize}
\item[(i)] $\fb$ has $n-2$ minimal generators of degree $a+n-3$ which have $n-3$
linear syzygies,
\item[(ii)] $\fb$ has a syzygy of degree $2a+d+n-5$.
\end{itemize}

Now we use a result of Rao (cf.\ \cite{Rao-res}, Theorem 2.5 and its
generalization
\cite{BR-sheaves}, Proposition 5.6, respectively). Combined with the
Cancellation principle (cf., e.g.,
\cite{Green-gin}, Corollary 1.21) we obtain that
$\TR_{n+1}(H^1_*(\cI_C), K)$ is a direct summand of $\TR_n (R/\gin I_C, K)$.
Therefore the
comparison with Lemma \ref{lem-res-stabil} shows that the syzygies of $\fb$ pointed
out in (i), (ii) above are the only syzygies of $\fb$. Hence $\fb$ has besides the
$n-2$ minimal generators of degree $a+n-3$ exactly one more generator $f$ which has
degree $a+d+n-5$.

Finally, the fact that the Hilbert function of $S/\fb$ has maximal growth in degree
$a+n-2$ implies that the $n-2$ minimal generators of degree $a+n-3$ of $\fb$ have a
common divisor, say $h$, of degree $a$. It follows that
$$
\fb = (h \cdot (x, y)^{n-3}, f).
$$
Since $S/\fb$ has finite length $f, h$ must be a regular sequence as claimed.
\end{proof}

\begin{remark} \label{rem-HR-res}
Although we excluded all extremal curves of degree $2$ and many of degree $3$ in
Theorem \ref{thm-HR-mod-extr} there is still some
information on the \HR\ of such curves.

(i)
Using the notation of the proof above its arguments  provide for extremal
curves with $d = 3, a \geq 1$ and $n \geq 4$ that $S/\fb$ must have a free
resolution of the form
$$
0 \to S^{n-3} (-a-n+2) \oplus S(-2a-n+2) \to S^{n-2} (-a-n+3) \oplus S(-a-n+2) \to S
\to S/\fb \to 0
$$
which is possibly not minimal. The latter case can really occur as Remark
\ref{rem-d3-more-HR} shows.

(ii)
A similar description is also possible for extremal curves $C \subset \PP^n$ of
degree $2$. Indeed, Proposition 3.1 in \cite{NNS1} implies that
then the \HR\ is
$$
H^1_* (\cI_C) \cong S/\fb (-g-1)
$$
where $S/\fb$ has a minimal free resolution of the form
$$
0 \to S^{n-3} (g-1) \oplus S(2g + n - 3) \to S^{n-1} (g)  \to S
\to S/\fb \to 0.
$$
\end{remark}
\medskip

The next goal is to show that extremal curves most often contain a
planar subcurve of largest possible degree.

\begin{proposition} \label{prop-subcurve}
Let $C \subset \PP^n$ be an extremal curve  satisfying either
\begin{itemize}
\item[(i)] $d \geq 5$ \quad or
\item[(ii)] $d = 4$ and $g < g_{max}$.
\end{itemize}
Then $C$ contains a planar subcurve of degree $d-1$.
\end{proposition}

\begin{proof}
In case (i) the claim follows by \cite{CGN3}, Corollary 4.4.
However, we will give a self-contained proof of our stronger
assertion.

Define the subscheme $Z \subset \PP^n$ as the intersection of the
quadrics containing the curve $C$. We may assume that $C$ and $Z$ are given
in generic coordinates. Then we claim that every minimal
basis of the ideal $I_Z$ is in fact a Gr\"obner basis.
Indeed, all the degree three S-polynomials of the generators of $I_Z$
reduce to zero because Proposition \ref{prop-gin-ex} shows that the
elements of the Gr\"obner basis of $I_C$ have either degree $2$ or
a degree $\geq \min \{d, d-1+a\} > 3$.

Using Proposition \ref{prop-gin-ex} again it follows that the
generic initial ideal of $Z$ is
$$
\gin I_Z = (x_0,\ldots,x_{n-4}) \cdot (x_0,\ldots,x_{n-1}) +
(x_{n-3}^2, x_{n-3} x_{n-2}).
$$
Computing the Hilbert polynomial of $R/\gin I_Z$ we see that $Z$
is a two-dimensional scheme of degree $1$, in other words $Z$ is
the union of a plane, $P$, and a scheme of dimension $\leq 1$. We
also may assume that $Z$ is in generic coordinates. Thus, we have
$$
\gin I_Z = (x_0,\ldots,x_{n-3}).
$$
Hence we obtain
$$
\iin (I_C \cap I_P) \subset \iin I_C \cap \iin I_P = \gin I_C \cap
\gin I_P  = \gin I_Z = \iin I_Z.
$$
Since we clearly have $I_Z \subset I_C \cap I_P$ it follows
$$
I_Z = I_C \cap I_P.
$$
Using the exact sequence
$$
0 \to R/I_Z \to R/I_C \oplus R/I_P \to R/(I_C + I_P) \to 0
$$
we see that $I_C + I_P$ defines a one-dimensional scheme which is
the union of a planar curve of degree $d-1$ and possibly a
zero-dimensional scheme. Our assertion follows.
\end{proof}

\begin{remark} Ellia  \cite{Ellia-curves} has shown that a
  non-degenerate curve $C \subset \PP^3$ of degree $d \geq 5$ is
  extremal if and only if $C$ contains a planar subcurve of degree
  $d-1$.
However, the converse of Proposition \ref{prop-subcurve}  is not true
  if  $n \geq 4$.
There are non-degenerate curves in $\PP^n, n \geq 4,$ of arbitrarily
  large degree  $d$ which do contain
a planar subcurve of degree $d-1$, but whose cohomology is smaller
than the one of extremal curves (cf.\  Example \ref{ex-non-max-Rao}).
\end{remark}

%%%%%%%%%%%%%%%%%%%%%%%%%%%%%%%%%%%%%%%%%%%%%%%%%%%%%%%%%%%%%

\section{Betti numbers of extremal curves}
\label{sec-Betti}

Lemma \ref{lem-res-stabil} provides upper bounds for the graded Betti
numbers of an  extremal curve $C$. The next result shows that
actually these bounds are attained in most cases.
The key will be the existence of another suitable subcurve of $C$.

\begin{proposition} \label{prop-Betti-numb}
Adopt the notation of Lemma \ref{lem-res-stabil}.
Let $C \subset \PP^n$ be an extremal curve  satisfying either
\begin{itemize}
\item[(i)] $d \geq 5$ \quad or
\item[(ii)] $d = 4$ and $g < g_{max}$.
\end{itemize}
Then  the minimal free resolution of $C$ has
the form
$$
0 \to R^{n-3} (-n-1) \oplus R(-d-a-n+2) \to F_{n-1} \to \ldots \to F_1 \to I
\to 0
$$
where  in case $a = g_{max} - g > 0$
$$
F_i = R^{\al_i}(-i-1) \oplus R^{\be_i}(-i-d+1) \oplus R^{\ga_i}(-i-d-a+2)
\quad (1 \leq i \leq n-1)
$$
and if  $a = 0$
$$
F_i = R^{\al_i}(-i-1) \oplus  R^{\be_i}(-i-d+2) \quad (1
\leq i \leq n-1).
$$
\end{proposition}

\begin{proof}
The Cancellation principle combined with Lemma
\ref{lem-res-stabil} shows that $C$ has a free resolution as
claimed. We have to show its minimality. This is obvious if $a =
0$ because then our assumption implies $d \geq 5$. Thus it remains
to consider the case $a \geq 1$.

Assuming that $C$ is in generic coordinates Proposition
\ref{prop-gin-ex} shows that the Gr\"obner basis of $I_C$ consists of
quadrics and two forms $g_1, g_2$ with
$$
\iin g_1 = x_{n-2}^d \quad \mbox{and} \quad \iin g_2 =
x_{n-2}^{d-1} x_{n-1}^a.
$$
Using the notation of the proof of Proposition \ref{prop-subcurve}
define the ideal $J \subset R$ by
$$
J := I_Z + g_1 R = [I_C]_2 R + g_1 R.
$$
Since $I_C \cap I_P = I_Z$ we obtain
$$
J \cap I_P = I_Z \subset I_P.
$$
It follows that $J + I_P = I_P + g_1 R$ defines a planar curve
of degree $d$.
Using the exact sequence
$$
0 \to R/I_Z \to R/J \oplus R/I_P \to R/(I_P + g_1 R) \to 0
$$
we see that $J$ has the same Hilbert function as $\gin (I_Z) +
x_{n-2}^d R$. Since obviously $\gin (I_Z) + x_{n-2}^d R \subset \iin
J$ we conclude
$$
\iin J = \gin (I_Z) + x_{n-2}^d R.
$$
Now $\iin (J)$ has a minimal free resolution of the form
$$
0 \to R^{n-3} (-n-1)  \to F_{n-1} \to \ldots \to F_1 \to I \to 0
$$
where,  using the notation of Lemma \ref{lem-res-stabil},
$$
F_i = R^{\al_i}(-i-1) \oplus R^{\be_i}(-i-d+1)  \quad (1 \leq i \leq
n-1).
$$
Since $d \geq 4$ any resolution of this shape must be minimal.
Thus, the graded Betti numbers of $J$ and $\iin (J)$ agree.
Observing that the minimal free resolution of $J$ injects into the
one of $I_C$ we see that the resolution described in the statement
 is minimal as claimed.
\end{proof}

Recall that an ideal $I \subset R$ is said to be {\it componentwise
linear} if the ideal $[I]_j R$ has a linear free resolution for
all integers $j$. Our last result implies.

\begin{corollary} \label{cor-linearly}
The homogeneous ideal of an extremal curve satisfying either
\begin{itemize}
\item[(i)] $d \geq 5$ \quad or
\item[(ii)] $d = 4$ and $g < g_{max}$.
\end{itemize}
is componentwise linear.
\end{corollary}

\begin{proof}
Let $C$ be an extremal curve satisfying our assumption.
Combining Proposition \ref{prop-gin-ex}, Lemma
\ref{lem-res-stabil} and Proposition \ref{prop-Betti-numb} we see
that the graded Betti numbers of $I_C$ and $\gin I_C$ agree. Hence
the main result of \cite{AHH} implies our claim.
\end{proof}

Finally, we show that Proposition \ref{prop-gin-ex}  can not improved
in general.

\begin{remark} \label{rem-gins}
We want to show that in case $d = 3, n \geq 4, a \geq 2$ the two
possiblities for the generic initial ideal of an extremal curve
allowed by Proposition \ref{prop-gin-ex} do really occur. Indeed, the
generic
initial ideal of the corresponding curve in Example \ref{ex-equations}
is
$$
(x_0,\ldots,x_{n-4}) \cdot (x_0,\ldots,x_{n-1}) + (x_{n-3}^2, x_{n-3}
x_{n-2}, x_{n-2}^3, x_{n-2}^{2} x_{n-1}^a).
$$
This follows by comparing the graded Betti numbers of the curve and
the ideals in Lemma \ref{lem-res-stabil} using the same arguments as
in the proof of Theorem \ref{thm-HR-mod-extr}.

Now consider the curve $C$ defined by the ideal
$$
I_C = (x_2,\ldots,x_n)^2 + (x_0^{a+1} x_3 + x_1^{a+1} x_4, x_0 x_i +
x_1 x_{i+1} \s 4 \leq i \leq n-1).
$$
Using Construction \ref{constr} we see that $C$ is an extremal curve
of degree $3$ and genus $g_{max} - a$. The ideal $I_C$ is
componentwise linear and
$$
\gin I_C = (x_0,\ldots,x_{n-5}) \cdot x_{n-1} +  (x_0,\ldots,x_{n-2})^2 +
(x_{n-4} x_{n-1}^{a+1}).
$$
\end{remark}
\smallskip 

The methods developed in this paper are rather general. Thus, we
believe that they can be succesfully applied to other related
problems. 

%%%%%%%%%%%%%%%%%%%%%%%%%%%%%%%%%%%%%%%%%%%%%%%%%%%%%%%%%%%%%%%%%

\end{document}